\documentclass[a4paper,12pt]{article}
\usepackage{amssymb}
\usepackage[latin1]{inputenc}     

\newcommand{\beq}{\begin{equation} }
\newcommand{\eqq}{\end{equation} }
\newcommand{\cuad}{{\sqcap\kern-.68em\sqcup}}

\newtheorem{remark}{Remark}[section]
\newcommand{\bremark}{\begin{remark} \em}
\newcommand{\eremark}{\end{remark} }

\def\beeq{\begin{equation}}
\def\eeq{\end{equation}}
\newcommand{\begeqaet}{\begin{eqnarray*}}
\newcommand{\eneqaet}{\end{eqnarray*}}

\let\Section=\section
\def\section{\setcounter{equation}{0}\Section}
\newtheorem{Lem}{Lemma}[section]
\newtheorem{Thm}{Theorem}[section]

\newtheorem{Remark}{Remark}[section]

\begin{document}
\begin{center}{\bf\Large Existence of solution for  perturbed fractional Hamiltonian systems}\medskip

\bigskip

\bigskip

{C\'esar Torres}

 Departamento de Ingenier\'{\i}a  Matem\'atica and
Centro de Modelamiento Matem\'atico,
 Universidad de Chile\\
 Santiago, Chile.\\
 {\sl  (ctorres@dim.uchile.cl)}

\end{center}

\medskip

\medskip
\medskip
\medskip
\medskip

\begin{abstract}
In this work we prove the existence of solution for a class of perturbed fractional Hamiltonian systems given by
\begin{eqnarray}\label{eq00}
-{_{t}}D_{\infty}^{\alpha}(_{-\infty}D_{t}^{\alpha}u(t)) - L(t)u(t)  + \nabla W(t,u(t)) = f(t),
\end{eqnarray}
where $\alpha \in (1/2, 1)$, $t\in \mathbb{R}$, $u\in \mathbb{R}^{n}$, $L\in C(\mathbb{R}, \mathbb{R}^{n^{2}})$ is a symmetric and positive definite matrix for all $t\in \mathbb{R}$, $W\in C^{1}(\mathbb{R}\times \mathbb{R}^{n}, \mathbb{R})$ and $\nabla W$ is the gradient of $W$ at $u$. The novelty of this paper is that, assuming $L$ is coercive at infinity we show that (\ref{eq00}) at least has one nontrivial solution.
\end{abstract}
\date{}

\setcounter{equation}{0}
\section{ Introduction}

In this paper, we shall be concerned with the existence of two solutions for the perturbed fractional Hamiltonian systems 
\begin{equation}\label{Eq00}
-{_{t}}D_{\infty}^{\alpha}(_{-\infty}D_{t}^{\alpha}u(t)) - L(t)u(t)  + \nabla W(t,u(t)) = f(t)
\end{equation}
where $t\in \mathbb{R}$ and $u\in \mathbb{R}^{n}$.

Fractional differential equations both ordinary and partial ones are applied in mathematical modeling of processes in
physics, mechanics, control theory, biochemistry, bioengineering and economics. Therefore the theory of fractional
differential equations is an area intensively developed during last decades \cite{OAJTMJS}, \cite{RH}, \cite{RMJK}, \cite{JSOAJTM}, \cite{BWMBPG}. The monographs \cite{AKHSJT}, \cite{KMBR}, \cite{IP}, enclose a review of methods of solving which are an extension of procedures from differential equations theory.

Recently, also equations including both - left and right fractional derivatives, are discussed. Let us point out that according to integration by parts formulas in fractional calculus , we obtain equations mixing left and right operators. Apart from their possible applications, equations with left and right derivatives are an interesting and new field in fractional differential equations theory. Some works in this topic can be founded in papers \cite{TABS}, \cite{DBJT}, \cite{MK} and its reference.

If $f = 0$ in (\ref{Eq00}), in \cite{CT} the author considered the following fractional Hamiltonian systems
\begin{eqnarray}\label{Eq02}
_{t}D_{\infty}^{\alpha}(_{-\infty}D_{t}^{\alpha}u(t)) + L(t)u(t)  = & \nabla W(t,u(t))
\end{eqnarray}
where $\alpha \in (1/2,1)$, $t\in \mathbb{R}$, $u\in \mathbb{R}^{n}$, $L\in C(\mathbb{R}, \mathbb{R}^{n\times n})$ is a symmetric matrix valued function for all $t\in \mathbb{R}$, $W\in C^{1}(\mathbb{R}\times \mathbb{R}^{n}, \mathbb{R})$ and $\nabla W(t, u(t))$ is the gradient of $W$ at $u$. Assuming that $L$ and $W$ satisfy the following hypotheses, I showed that (\ref{Eq02}) possesses at least one nontrivial solution via Mountain pass theorem. Explicitly,
\begin{itemize}
\item[$(L)$] $L(t)$ is positive definite symmetric matrix for all $t\in \mathbb{R}$ and there exists an $l\in C(\mathbb{R}, (0,\infty))$ such that $l(t) \to +\infty$ as $t \to \infty$ and
    \begin{equation}\label{Eq03}
    (L(t)x,x) \geq l(t)|x|^{2},\;\;\mbox{for all}\;t\in \mathbb{R}\;\;\mbox{and}\;\;x\in \mathbb{R}^{n}.
    \end{equation}
\item[$(W_{1})$] $W\in C^{1}(\mathbb{R} \times \mathbb{R}^{n}, \mathbb{R})$ and there is a constant $\mu >2$ such that
$$
0< \mu W(t,x) \leq (x, \nabla W(t,x)),\;\;\mbox{for all}\;t\in \mathbb{R}\;\;\mbox{and}\;x\in\mathbb{R}^{n}\setminus \{0\}.
$$
\item[$(W_{2})$] $|\nabla W(t,x)| = o(|x|)$ as $x\to 0$ uniformly with respect to $t\in \mathbb{R}$.
\item[$(W_{3})$] There exists $\overline{W} \in C(\mathbb{R}^{n}, \mathbb{R})$ such that
$$
|W(t,x)| + |\nabla W(t,x)| \leq |\overline{W(x)}|\;\;\mbox{for every}\;x\in \mathbb{R}^{n}\;\mbox{and}\;t\in \mathbb{R}.
$$
\end{itemize}


\noindent
Very recently Zhang and Yuan \cite{ZZRY}, using the genus properties of critical point theory, the authors generalized the result of \cite{CT} and established some new criterion to guarantee the existence of infinitely many solutions of (\ref{Eq02}) for the case that $W(t,u)$ is subquadratic as $|u| \to +\infty$. Explicitly, $L$ satisfies $(L)$ and the potential $W(t,u)$ is supposed to satisfy the following conditions:
\begin{itemize}
\item[$(HS)_{1}$] $W(t,0) = 0$ for all $t\in \mathbb{R}$, $W(t,u) \geq a(t)|u|^{\theta}$ and $|\delta W(t,u)| \leq b(t)|u|^{\theta -1}$ for all $(t,u)\in \mathbb{R} \times \mathbb{R}^{n}$, where $ \theta < 2$ is a constant, $a :\mathbb{R} \to \mathbb{R}^{+}$ is a bounded continuous function and $b:\mathbb{R} \to \mathbb{R}^{+}$ is a continuous function such that $b\in L^{\frac{2}{2-\theta}}(\mathbb{R})$;
\item[$(HS)_{2}$] There is a constant $1 < \sigma \leq \theta < 2$ such that
$$
(W(t,u), u) \leq \sigma W(t,u)\quad \mbox{for all t}\in \mathbb{R}\;\;\mbox{and}\;\;u\in \mathbb{R}^{n}\setminus \{0\};
$$
\item[$(HS)_{3}$] $W(t,u)$ is even in $u$, i.e. $W(t,u) = W(t,-u)$ for all $t\in \mathbb{R}$ and $u\in \mathbb{R}^{n}$.
\end{itemize}

In this paper we extend the result of Torres \cite{CT} and Zhang and Yuan \cite{ZZRY} to the case $f \neq 0$. For that purpose throughout the paper, $(.,.): \mathbb{R}^{n} \times \mathbb{R}^{n} \to \mathbb{R}$ denotes the standar scalar product in $\mathbb{R}^{n}$ and $|.|$ is the induced norm. Let
$$
X^{\alpha} = \left\{ u\in H^{\alpha}(\mathbb{R}, \mathbb{R}^{n}):\;\;\int_{\mathbb{R}} [|_{-\infty}D_{t}^{\alpha}u(t)|^{2} + (L(t)u(t), u(t))]dt < \infty\right\}.
$$
$X^{\alpha}$ is a Hilbert space under the norm
$$
\|u\|_{X^{\alpha}} = \left( \int_{\mathbb{R}} [|_{-\infty}D_{t}^{\alpha}u(t)|^{2} + (L(t)u(t), u(t))]dt \right)^{1/2}.
$$

\noindent 
From ($L$) it follows that there is a constant $C_{e} >0$ such that for every $u\in X^{\alpha}$,
\begin{equation}\label{emb}
\|u\|_{H^{\alpha}} \leq C_{\epsilon}\|u\|_{X^{\alpha}},
\end{equation}
see \cite{CT} for more details. Set 
$$
M = \max_{|u|=1} \overline{W}(u).
$$
We will also assume that
\begin{itemize}
\item[($W_{f}$)] $M< \frac{1}{2C_{\alpha}^{2}C_{e}^{2}}$ and $f:\mathbb{R} \to \mathbb{R}^{n}$ is a continuous square integrable function such that 
$$
\|f\|_{L^{2}} < \frac{1}{2C_{\alpha}^{2}C_{e}^{2}} - M,
$$
where $C_{\alpha}>0$ (see Theorem \ref{FDEtm01} below). 
\end{itemize}

We are going to prove the following theorem.
\begin{Thm}\label{tm01}
Suppose that $(L)-(W_{3})$ and ($W_{f}$) hold, then (\ref{Eq00}) possesses at least two nontrivial solution in $X^{\alpha}$.
\end{Thm}

Under this assumption, the problem of existence of solutions is much more delicate, because extra difficulties arise in studying the properties of the corresponding action functional $I : X^{\alpha} \to  R$.

The problem here is as follows. We are given two sequences of almost critical points in $X^{\alpha}$. The first one, obtained by Ekeland's variational principle, is contained in a small ball centered at $0$. Using the mountain pass geometry of the action functional, the existence of the second sequence is established. Both sequences are weakly convergent in $X^{\alpha}$. The question is whether their weak limits are equal to each other or they define two geometrically distinct solutions of (\ref{Eq00}). The PS-condition is enough to obtain two solutions.

The assumption $(L)$ ensures the PS-condition at each level. In fact one needs the PS-condition only at two levels and therefore it is tempting to seek for weaker compactness assumptions.

The rest of the paper is organized as follows: in section 2, subsection 2.1, we describe the Liouville-Weyl fractional calculus; in subsection 2.2 we introduce the fractional space that we use in our work and some proposition are proven which will aid in our analysis. In section 3, we will prove theorem \ref{tm01}.

\section{Preliminary Results}

\subsection{Liouville-Weyl Fractional Calculus}

The Liouville-Weyl fractional integrals of order $0<\alpha < 1$ are defined as
\begin{equation}\label{LWeq01}
_{-\infty}I_{x}^{\alpha}u(x) = \frac{1}{\Gamma (\alpha)} \int_{-\infty}^{x}(x-\xi)^{\alpha - 1}u(\xi)d\xi
\end{equation}
\begin{equation}\label{LWeq02}
_{x}I_{\infty}^{\alpha}u(x) = \frac{1}{\Gamma (\alpha)} \int_{x}^{\infty}(\xi - x)^{\alpha - 1}u(\xi)d\xi
\end{equation}
The Liouville-Weyl fractional derivative of order $0<\alpha <1$ are defined as the left-inverse operators of the corresponding Liouville-Weyl fractional integrals
\begin{equation}\label{LWeq03}
_{-\infty}D_{x}^{\alpha}u(x) = \frac{d }{d x} {_{-\infty}}I_{x}^{1-\alpha}u(x)
\end{equation}
\begin{equation}\label{LWeq04}
_{x}D_{\infty}^{\alpha}u(x) = -\frac{d }{d x} {_{x}}I_{\infty}^{1-\alpha}u(x)
\end{equation}
The definitions (\ref{LWeq03}) and (\ref{LWeq04}) may be written in an alternative form:
\begin{equation}\label{LWeq05}
_{-\infty}D_{x}^{\alpha}u(x) = \frac{\alpha}{\Gamma (1-\alpha)} \int_{0}^{\infty}\frac{u(x) - u(x-\xi)}{\xi^{\alpha + 1}}d\xi
\end{equation}
\begin{equation}\label{LWeq05}
_{x}D_{\infty}^{\alpha}u(x) = \frac{\alpha}{\Gamma (1-\alpha)} \int_{0}^{\infty}\frac{u(x) - u(x+\xi)}{\xi^{\alpha + 1}}d\xi
\end{equation}

\noindent
We establish the Fourier transform properties of the fractional integral and fractional differential operators. Recall that the Fourier transform $\widehat{u}(w)$ of $u(x)$ is defined by
$$
\widehat{u}(w) = \int_{-\infty}^{\infty} e^{-ix.w}u(x)dx.
$$
Let $u(x)$ be defined on $(-\infty, \infty)$. Then the Fourier transform of the Liouville-Weyl integral and differential operator satisfies
\begin{equation}\label{LWeq06}
\widehat{ _{-\infty}I_{x}^{\alpha}u(x)}(w) = (iw)^{-\alpha}\widehat{u}(w)
\end{equation}
\begin{equation}\label{LWeq07}
\widehat{ _{x}I_{\infty}^{\alpha}u(x)}(w) = (-iw)^{-\alpha}\widehat{u}(w)
\end{equation}
\begin{equation}\label{LWeq08}
\widehat{ _{-\infty}D_{x}^{\alpha}u(x)}(w) = (iw)^{\alpha}\widehat{u}(w)
\end{equation}
\begin{equation}\label{LWeq09}
\widehat{ _{x}D_{\infty}^{\alpha}u(x)}(w) = (-iw)^{\alpha}\widehat{u}(w)
\end{equation}
\subsection{Fractional Derivative Spaces}

In this section we introduce some fractional spaces for more detail see \cite{VEJR}.

\noindent
Let $\alpha > 0$. Define the semi-norm
$$
|u|_{I_{-\infty}^{\alpha}} = \|_{-\infty}D_{x}^{\alpha}u\|_{L^{2}}
$$
and norm
\begin{equation}\label{FDEeq01}
\|u\|_{I_{-\infty}^{\alpha}} = \left( \|u\|_{L^{2}}^{2} + |u|_{I_{-\infty}^{\alpha}}^{2} \right)^{1/2},
\end{equation}
and let
$$
I_{-\infty}^{\alpha} (\mathbb{R}, \mathbb{R}^{n}) = \overline{C_{0}^{\infty}(\mathbb{R}, \mathbb{R}^{n})}^{\|.\|_{I_{-\infty}^{\alpha}}}.
$$
Now we define the fractional Sobolev space $H^{\alpha}(\mathbb{R}, \mathbb{R}^{n})$ in terms of the Fourier transform. Let $0< \alpha < 1$, let the semi-norm
\begin{equation}\label{FDEeq02}
|u|_{\alpha} = \||w|^{\alpha}\widehat{u}\|_{L^{2}}
\end{equation}
and norm
$$
\|u\|_{\alpha} = \left( \|u\|_{L^{2}}^{2} + |u|_{\alpha}^{2} \right)^{1/2},
$$
and let
$$
H^{\alpha}(\mathbb{R}, \mathbb{R}^{n}) = \overline{C_{0}^{\infty}(\mathbb{R}, \mathbb{R}^{n})}^{\|.\|_{\alpha}}.
$$

\noindent
We note a function $u\in L^{2}(\mathbb{R}, \mathbb{R}^{n})$ belong to $I_{-\infty}^{\alpha}(\mathbb{R}, \mathbb{R}^{n})$ if and only if
\begin{equation}\label{FDEeq03}
|w|^{\alpha}\widehat{u} \in L^{2}(\mathbb{R}, \mathbb{R}^{n}).
\end{equation}
Especially
\begin{equation}\label{FDEeq04}
|u|_{I_{-\infty}^{\alpha}} = \||w|^{\alpha}\widehat{u}\|_{L^{2}}.
\end{equation}
Therefore $I_{-\infty}^{\alpha}(\mathbb{R}, \mathbb{R}^{n})$ and $H^{\alpha}(\mathbb{R}, \mathbb{R}^{n})$ are equivalent with equivalent semi-norm and norm. Analogous to $I_{-\infty}^{\alpha}(\mathbb{R}, \mathbb{R}^{n})$ we introduce $I_{\infty}^{\alpha}(\mathbb{R}, \mathbb{R}^{n})$. Let the semi-norm
$$
|u|_{I_{\infty}^{\alpha}} = \|_{x}D_{\infty}^{\alpha}u\|_{L^{2}}
$$
and norm
\begin{equation}\label{FDEeq05}
\|u\|_{I_{\infty}^{\alpha}} = \left( \|u\|_{L^{2}}^{2} + |u|_{I_{\infty}^{\alpha}}^{2} \right)^{1/2},
\end{equation}
and let
$$
I_{\infty}^{\alpha}(\mathbb{R}, \mathbb{R}^{n}) = \overline{C_{0}^{\infty}(\mathbb{R}, \mathbb{R}^{n})}^{\|.\|_{I_{\infty}^{\alpha}}}.
$$
Moreover $I_{-\infty}^{\alpha}(\mathbb{R}, \mathbb{R}^{n})$ and $I_{\infty}^{\alpha}(\mathbb{R}, \mathbb{R}^{n})$ are equivalent , with equivalent semi-norm and norm \cite{VEJR}.

Now we give the prove of the Sobolev lemma.
\begin{Thm}\label{FDEtm01}
\cite{} If $\alpha > \frac{1}{2}$, then $H^{\alpha}(\mathbb{R}, \mathbb{R}^{n}) \subset C(\mathbb{R}, \mathbb{R}^{n})$ and there is a constant $C=C_{\alpha}$ such that
\begin{equation}\label{FDEeq06}
\sup_{x\in \mathbb{R}} |u(x)| \leq C \|u\|_{\alpha}
\end{equation}
\end{Thm}
\begin{Remark}\label{FDEnta01}
If $u\in H^{\alpha}(\mathbb{R}, \mathbb{R}^{n})$, then $u\in L^{q}(\mathbb{R}, \mathbb{R}^{n})$ for all $q\in [2,\infty]$, since
$$
\int_{\mathbb{R}} |u(x)|^{q}dx \leq \|u\|_{\infty}^{q-2}\|u\|_{L^{2}}^{2}
$$
\end{Remark}

\noindent
We introduce a new fractional spaces. Let
$$
X^{\alpha} = \left\{ u\in H^{\alpha}(\mathbb{R}, \mathbb{R}^{n})|\;\;\int_{\mathbb{R}} |_{-\infty}D_{t}^{\alpha}u(t)|^{2} + (L(t)u(t), u(t)) dt < \infty  \right\}
$$
The space $X^{\alpha}$ is a Hilbert space with the inner product
$$
\langle u,v \rangle_{X^{\alpha}} = \int_{\mathbb{R}} (_{-\infty}D_{t}^{\alpha}u(t) , \; _{-\infty}D_{t}^{\alpha}v(t)) + (L(t)u(t), v(t))dt
$$
and the corresponding norm
$$
\|u\|_{X^{\alpha}}^{2} = \langle u,u \rangle_{X^{\alpha}}
$$
\begin{Lem}\label{FDElm01}
\cite{CT} Suppose $L$ satisfies ($L$). Then $X^{\alpha}$ is continuously embedded in $H^{\alpha}(\mathbb{R},\mathbb{R}^{n})$.
\end{Lem}

\begin{Lem}\label{FDElm02}
\cite{CT} Suppose $L$ satisfies ($L$). Then the imbedding of $X^{\alpha}$ in $L^{2}(\mathbb{R}, \mathbb{R}^{n})$ is compact.
\end{Lem}

\begin{Lem}\label{FDElm03}
\cite{CT} 
\begin{equation}\label{FDEeq07}
W(t,u) \geq W(t, \frac{u}{|u|})|u|^{\mu},\;\;|u|\geq 1
\end{equation}
and
\begin{equation}\label{FDEeq08}
W(t,u) \leq W(t, \frac{u}{|u|})|u|^{\mu},\;\;|u|\leq 1
\end{equation}
\end{Lem}

\begin{Remark}\label{FDE-nta01}

By Lemma \ref{FDElm03}, we have
\begin{equation}\label{FDEeq09}
W(t, u) = o(|u|^{2})\;\mbox{as}\;u\to 0\;\mbox{uniformly in}\;t\in\mathbb{R}
\end{equation}
In addition, by $(W_{2})$, we have, for any $u\in \mathbb{R}^{n}$ such that $|u| \leq M_{1}$, there exists some constant $d>0$ (dependent on $M_{1}$) such that
\begin{equation}\label{FDEeq10}
|\nabla W(t,u(t))| \leq d|u(t)|
\end{equation}
\end{Remark}

\noindent

\begin{Lem}\label{FDElm04}
\cite{CT} Suppose that (L), ($W_{1}$)-($W_{2}$) are satisfied. If $u_{k} \rightharpoonup u$ in $X^{\alpha}$, then $\nabla W(t, u_{k}) \to \nabla W(t, u)$ in $L^{2}(\mathbb{R}, \mathbb{R}^{n})$.
\end{Lem}


Now we introduce more notations and some necessary definitions. Let $\mathfrak{B}$ be a real Banach space, $I\in C^{1}(\mathfrak{B},\mathbb{R})$, which means that $I$ is a continuously Fréchet-differentiable functional defined on $\mathfrak{B}$. Recall that $I\in C^{1}(\mathfrak{B},\mathbb{R})$ is said to satisfy the (PS) condition if any sequence $\{u_{k}\}_{k\in \mathbb{N}} \in \mathfrak{B}$, for which $\{I(u_{k})\}_{k\in \mathbb{N}}$ is bounded and $I'(u_{k}) \to 0$ as $k\to +\infty$, possesses a convergent subsequence in $\mathfrak{B}$.

Moreover, let $B_{r}$ be the open ball in $\mathfrak{B}$ with the radius $r$ and centered at $0$ and $\partial B_{r}$ denote its boundary. We obtain the existence of homoclinic solutions of (\ref{Eq02}) by use of the following well-known Mountain Pass Theorems, see \cite{PR}.
\begin{Thm}\label{FDEtm02}
Let $\mathfrak{B}$ be a real Banach space and $I\in C^{1}(\mathfrak{B}, \mathbb{R})$ satisfying (PS) condition. Suppose that $I(0) = 0$ and
\begin{itemize}
\item[i.] There are constants $\rho , \beta >0$ such that $I|_{\partial B_{\rho}} \geq \beta$, and
\item[ii.] There is and $e\in \mathfrak{B} \setminus \overline{B_{\rho}}$ such that $I(e)\leq 0$.
\end{itemize}
Then $I$ possesses a critical value $c\geq \beta$. Moreover $c$ can be characterized as
$$
c = \inf_{\gamma \in \Gamma} \max_{s\in [0,1]} I(\gamma (s))
$$
where
$$
\Gamma = \{\gamma \in C([0,1] , \mathfrak{B}):\;\;\gamma (0) = 0,\;\;\gamma (1) = e\}
$$
\end{Thm}
\section{Proof of Theorem \ref{tm01}}

Now we are going to establish the corresponding variational framework to obtain the existence of solutions for (\ref{Eq02}). Define the functional $I: X^{\alpha} \to \mathbb{R}$ by
\begin{eqnarray}\label{TMeq01}
I(u) & = & \int_{\mathbb{R}} \left[ \frac{1}{2}|_{-\infty}D_{t}^{\alpha}u(t)|^{2} + \frac{1}{2}(L(t)u(t),u(t)) - W(t,u(t)) + (f(t), u(t))\right]dt \nonumber\\
     & = & \frac{1}{2}\|u\|_{X^{\alpha}}^{2} - \int_{\mathbb{R}} W(t,u(t))dt+ \int_{\mathbb{R}} (f(t), u(t))dt.
\end{eqnarray}

\noindent
Then $I\in C^{1}(X^{\alpha}, \mathbb{R})$ and it is easy to check that
\begin{equation}\label{TMeq02}
I'(u)v = \int_{\mathbb{R}} \left[ (_{-\infty}D_{t}^{\alpha}u(t), _{-\infty}D_{t}^{\alpha}v(t)) + (L(t)u(t),v(t)) - (\nabla W(t,u(t)),v(t)) + (f(t), v(t))\right]dt
\end{equation}
for all $u,v \in X^{\alpha}$, which yields that
\begin{equation}\label{TMeq03}
I'(u)u = \|u\|_{X^{\alpha}}^{2} - \int_{\mathbb{R}}(\nabla W(t,u(t)), u(t))dt + \int_{\mathbb{R}}(f(t), u(t)).
\end{equation}

In order to prove Theorem \ref{tm01} we use the mountain pass theorem and Ekeland's variational principle \cite{JMMW}, \cite{PR}.  The proof will be divided into a sequence of lemmas. 
\begin{Lem}\label{TMlm01}
Suppose that $(L)-(W_{3})$ and $(W_{f})$ holds. Then $I$ satisfies the PS-condition. 
\end{Lem}

\noindent
{\bf Proof.} Assume that $(u_{k})_{k\in \mathbb{N}} \in X^{\alpha}$ is a sequence such that 
\begin{equation}\label{TMeq04}
I(u_{k}) \to c \quad \mbox{ and}\quad  I'(u_{k}) \to 0\;\;\mbox{ as}\;\; k \to +\infty
\end{equation}
We have to show that $\{u_{k}\}_{k\in \mathbb{N}}$ possesses a convergent subsequence in $X^{\alpha}$ 

\noindent
We firstly prove that $\{u_{k}\}_{k\in \mathbb{N}}$ is bounded in $X^{\alpha}$. By (\ref{TMeq01}), (\ref{TMeq03}) we get
\begin{eqnarray}\label{TMeq05}
I(u_{k})-\frac{1}{\mu}I'(u_{k})u_{k} & = & \left( \frac{1}{2}  - \frac{1}{\mu}\right)\|u_{k}\|_{X^{\alpha}}^{2} \nonumber\\
&& - \int_{\mathbb{R}} \left[W(t,u_{k}(t)) - \frac{1}{\mu}(\nabla W(t,u_{k}(t)), u_{k}(t))\right]dt \nonumber\\
& & + \left( 1- \frac{1}{\mu}\right) \int_{\mathbb{R}} (f(t), u_{k}(t))dt .
\end{eqnarray}
From ($W_{1}$) and Lemma \ref{FDElm01}, it follows that
$$
I(u_{k}) - \frac{1}{\mu} I'(u_{k})u_{k} \geq \left( \frac{1}{2} - \frac{1}{\mu} \right) \|u_{k}\|_{X^{\alpha}}^{2} -  \left( C_{e} - \frac{C_{e}}{\mu} \right) \|f\|_{L^{2}} \|u_{k}\|_{X^{\alpha}},
$$
where $C_{e}$ denote the constant of the continuous embedding given by Lemma \ref{FDElm01}. On the other hand, by (\ref{TMeq04}), there is $k_{0}\in \mathbb{N}$ such that for $k\geq k_{0}$,
$$
c + 1 + \|u_{k}\|_{X^{\alpha}} \geq I(u_{k}) - \frac{1}{\mu}I'(u_{k})u_{k}.
$$
In consequence, since $\mu >2$, $\{u_{k}\}_{k\in \mathbb{N}}$ is bounded in $X^{\alpha}$. Since $X^{\alpha}$ is a Hilbert space, passing to a subsequence if necessary, it can be assumed that $u_{k} \rightharpoonup u$ in $X^{\alpha}$ and hence, by Lemma \ref{FDElm02}, $u_{k} \to u$ in $L^{2}(\mathbb{R},\mathbb{R}^{n})$. It follows from the definition of $I$ that
\begin{eqnarray}\label{TMeq12}
\|u_{k} - u\|_{X^{\alpha}}^{2} & = & (I'(u_{k}) - I'(u))(u_{k} - u) \nonumber\\
&& + \int_{\mathbb{R}}[\nabla W(t,u_{k}) - \nabla W(t,u)](u_{k} - u)dt.
\end{eqnarray}
Since $u_{k} \to u$ in $L^{2}(\mathbb{R}, \mathbb{R}^{n})$, by Lemma \ref{FDElm03} $\nabla W(t, u_{k}(t)) \to \nabla W(t,u(t))$ in $L^{2}(\mathbb{R}, \mathbb{R}^{n})$. Hence
$$
\int_{\mathbb{R}} ( \nabla W(t,u_{k}(t))-\nabla W(t,u(t)), u_{k}(t)-u(t) )dt \to 0
$$
as $k\to +\infty$. So (\ref{TMeq12}) implies
$$
\|u_{k} - u\|_{X^{\alpha}} \to 0\;\;\mbox{as}\;\;k\to +\infty.
$$
$\Box$

\begin{Lem}\label{TMlm02}
Suppose that $(L)-(W_{3})$ and ($W_{f}$) holds. There are $\rho >0$ and $\beta$ such that $I(u) \geq \beta$ for $\|u\|_{X^{\alpha}} = \rho$
\end{Lem}

\noindent
{\bf Proof.} Let $\rho = \frac{1}{C_{\alpha} C_{e}}$. Assume that $u\in X^{\alpha}$ and $\|u\|_{\alpha} \leq \rho$. By Theorem \ref{FDEtm01} and (\ref{emb})
$$
\|u\|_{\infty} \leq C_{\alpha}C_{e} \|u\|_{X^{\alpha}} \leq 1
$$
By (\ref{FDEeq08}) and ($W_{3}$), we get
\begin{eqnarray}\label{TMeq07}
I(u) & = & \frac{1}{2}\|u\|_{X^{\alpha}}^{2} - \int_{\mathbb{R}} W(t,u(t))dt+ \int_{\mathbb{R}} (f(t), u(t))dt \nonumber\\
& \geq & \frac{1}{2}\|u\|_{X^{\alpha}}^{2} - \int_{\mathbb{R}} W(t, \frac{u}{|u|})|u|^{\mu}dt - \|f\|_{L^{2}}\|u\|_{L^{2}} \nonumber \\
& \geq & \frac{1}{2}\|u\|_{X^{\alpha}}^{2} - M\|u\|_{\infty}^{\mu - 2} \int_{\mathbb{R}}|u|^{2}dt - C_{\alpha}C_{e}\|f\|_{L^{2}}\|u\|_{X^{\alpha}} \nonumber\\
& \geq & \frac{1}{2}\|u\|_{X^{\alpha}}^{2} - MC_{\alpha}^{2}C_{e}^{2} \|u\|_{X^{\alpha}}^{2} - C_{\alpha}C_{e} \|f\|_{L^{2}} \|u\|_{X^{\alpha}}.
\end{eqnarray}
Hence, if $\|u\|_{X^{\alpha}} = \rho$, we have
$$
I(u) \geq \frac{1}{2C_{\alpha}^{2}C_{e}^{2}} - M - \|f\|_{L^{2}} = \beta >0,
$$
by ($W_{f}$). $\Box$


\begin{Lem}\label{TMlm03}
Suppose that $(L)-(W_{3})$ and ($W_{f}$) holds. There is $e\in X^{\alpha} \setminus \overline{B(0, \rho)}$ such that $I(e) \leq 0$, where $B(0, \rho)$ is a ball in $X^{\alpha}$ of radius $\rho$ centered at $0$ and $\rho$ is given by Lemma \ref{TMlm02}.
\end{Lem}

\noindent
{\bf Proof.} Fix $u\in X^{\alpha}$ such that $|u(t)| = 1$ for all $t\in [0,1]$ and assume that $\sigma \geq 1$. Then by (\ref{FDEeq07})
\begin{eqnarray*}
I(\sigma u) & = & \frac{\sigma^{2}}{2}\|u\|_{X^{\alpha}}^{2} - \int_{\mathbb{R}}W(t,\sigma u(t))dt + \sigma \int_{\mathbb{R}} (f(t), u(t))dt\\
& \leq & \frac{\sigma^{2}}{2}\|u\|_{X^{\alpha}}^{2} - \sigma^{\mu} \int_{0}^{1} W\left( t, \frac{u(t)}{|u(t)|}\right)|u(t)|^{\mu} dt + \sigma \int_{\mathbb{R}} (f(t), u(t))dt\\
& \leq & \frac{\sigma^{2}}{2}\|u\|_{X^{\alpha}}^{2} - \sigma^{\mu}m\int_{0}^{1}|u(t)|^{\mu}dt + \sigma \int_{\mathbb{R}} (f(t), u(t))dt
\end{eqnarray*}
$\Box$
where
$$
m = \min_{t\in [0,1], |u|=1} W(t,u).
$$
Since $\mu >2$, $I(\sigma u) \to -\infty$ as $\sigma \to \infty$. Hence there is $\sigma \geq 1$ such that $\|\sigma u\|_{X^{\alpha}} > \rho$ and $I(\sigma u) \leq 0$. $\Box$

\medskip
\noindent
{\bf Proof of Theorem \ref{tm01}} Since $I(0) = 0$ and $I$ satisfies Lemmas \ref{TMlm01} - \ref{TMlm03}, it follows by the mountain pass theorem that $I$ has a critical value $c$ given by
$$
c = \inf_{\gamma \in \Gamma} \max_{t\in [0,1]} I(\gamma (t)),
$$
where 
$$
\Gamma = \{\gamma \in C([0,1], X^{\alpha}):\;\;\gamma (0) = 0, \; I(\gamma (1)) < 0\}.
$$
By definition, it follows that $c \geq \beta >0$. By (\ref{TMeq07}), $I$ is bounded from below on $\overline{B(0, \rho)}$. Let 
\begin{equation}\label{TMeq08}
c_{1} = \inf_{\|u\|_{X^{\alpha}} \leq \rho} I(u).
\end{equation}
Since $I(0) = 0$, $c_{1} \leq 0$. Thus $c_{1} \leq c$. By Ekeland's variational principle, there is a minimizing sequence $\{w_{k}\}_{k\in \mathbb{N}} \subset \overline{B(0,\rho)}$ such that 
$$
I(w_{k}) \to c_{1}\quad\mbox{and}\quad I'(w_{k})\to 0,
$$
as $k\to \infty$. From Lemma \ref{TMlm01}, $c_{1}$ is a critical value of $I$. Consequently, $I$ has at least two critical points. $\Box$


\noindent {\bf Acknowledgements:}
C.T. was  partially supported by MECESUP 0607 and CMM.

\end{document}